\documentclass[12pt]{article}
\usepackage{amsmath,amsfonts,amsthm,amscd}
\newtheorem{theorem}{Theorem}

\newtheorem{proposition}{Proposition}
\newtheorem{lemma}{Lemma}
\newtheorem{definition}{Definition}
\newtheorem{corollary}{Corollary}

\numberwithin{equation}{section}
\numberwithin{lemma}{section}
\numberwithin{claim}{section}
\input{diagrams.tex}
\newcommand{\bull}{\ensuremath{{}\bullet{}}}
\newcommand{\gc}{\ensuremath{SL(N+1,\mathbb{C})}}

\newcommand{\lam}{\ensuremath{\lambda}}

\newcommand{\gr}{\ensuremath{\mathbb{G}(N-n,\mathbb{C}^{N+1})}}

\newcommand{\cpn}{\ensuremath{\mathbb{P}^{N}}}
\newcommand{\dlb}{\ensuremath{\overline{\partial}}}
\newcommand{\dl}{\ensuremath{\partial}}
\newcommand{\ra}{\ensuremath{\longrightarrow}}
\newcommand{\ba}{\ensuremath{\begin{align*}}}
\newcommand{\ea}{\ensuremath{\end{align*}}}
\newcommand{\vp}{\ensuremath{\varphi}}
\newcommand{\vpt}{\ensuremath{\varphi_{t}}}
\newcommand{\om}{\ensuremath{\omega}}
\newcommand{\tp}{\ensuremath{\frac{\sqrt{-1}}{2\pi}}}
 \newarrow{to}---->
\newarrow{mto}|--->
\newarrowtail{<=}\Rightarrow\Leftarrow\Uparrow\Downarrow
\newarrow{bboth}{<=}==={=>}
\newarrow{inject}{hooka}---{vee}
\newarrow{surject}----{>>}

\begin{document}
\title{\sc{Higher Energies in K\"ahler Geometry I}}
\author{\sc{Sean Timothy Paul}\thanks {The University of Wisconsin Madison. The author is supported by an NSF DMS grant  0505059.}}
 
\date{ } 
\maketitle
\vspace{-5mm}
 
 {\begin{abstract}{ Let $X\hookrightarrow \cpn $ be a smooth  complex projective variety of dimension $n$. Let $\lambda$ be an algebraic one parameter subgroup of $G:=\gc$. 
Let $ 0\leq l\leq n+1$. We associate to the coefficients $F_{l}(\lambda)$ of the normalized weight of $\lambda$ on the $mth$ Hilbert point of $X$  new energies $F_{\om,l}(\vp)$.
The (logarithmic) asymptotics of  $F_{\om,l}(\vp)$ along the potential deduced from $\lambda$ is the weight  $F_{l}(\lambda)$.  $F_{\om,l}(\vp)$ reduces to the Aubin energy when $l=0$ and the K-Energy map of Mabuchi when $l=1$. When $l\geq 2$ $F_{\om,l}(\vp)$ coincides (modulo lower order terms) with  the  functional  $E_{\om,l-1}(\vp)$ introduced by X.X. Chen and G.Tian.}
\end{abstract}   

\begin{center}{\S 0 \sc{The Standard Energies of K\"ahler Geometry}}\end{center} 
  Recall that Mabuchi's {K-energy map} (see \cite{mabuchi} ) is given by
\begin{align*}
    \nu_{\omega}(\varphi):= -\frac{1}{V}\int_{0}^{1}\int_{X}\dot{\varphi_{t}}(\mbox{Scal}(\varphi_{t})-\mu)\omega_{t}^{n}dt \ .
\end{align*}
 $\vp_{t}$ denotes a path in $P(X,\om)$, $\mbox{Scal}(\varphi_{t})$ denotes the scalar curvature of the metric $\omega +\sqrt{-1}\dl\dlb\varphi_{t}$ and $\mu$ denotes the average of the scalar curvature. Critical points of the K-Energy are metrics of constant scalar curvature.  In \cite{bandmab} Bando and Mabuchi have proved that the K-Energy is \emph{bounded from below} provided that  $(X,\om)$ admits a K\"ahler Einstein metric (in this case it is required that $\om=c_{1}(K^{-1}_{X})$). This is noteworthy as the K Energy map is (essentially) the \emph{difference} of two \emph{positive} terms as follows.
 \begin{align*}
 &\nu_{\omega}(\varphi)=  
\frac{1}{V}\int_{M}\log\left(
\frac{\omega_{\varphi}^{n}}{\omega^{n}}\right)\omega_{\varphi}^{n}
 -\frac{\mu}{n}\left(I_{\omega}(\varphi)-J_{\omega}(\varphi)\right)
+\frac{1}{V}\int_{M}h_{\om}\om_{\vp}^{n}\\
 & I_{\omega}(\varphi):= \frac{1}{V}\int_{X}\varphi(\omega^{n}-\omega_{\varphi}^{n})\\
 & J_{\omega}(\varphi):= \frac{\sqrt{-1}}{V}\int_{X}\sum_{i=0}^{n-1}\frac{i+1}{n+1}\dl\varphi \wedge \dlb
\varphi\wedge \omega^{i}\wedge \om_{\vp}^{n-i-1}.\\
\end{align*}
The work of Bando and Mabuchi has been extended to \emph{any} K\"ahler class by X.X. Chen and G. Tian.

Recall that $\nu_{\omega}$ is \emph{proper}  (see \cite{psc}) provided there is an increasing function $f:\mathbb{R}\rightarrow \mathbb{R}$ such that $\nu_{\omega}(\varphi)\geq f(J_{\omega}(\varphi))$ for all $\varphi\in P(X,\omega)$.  This notion is due to Gang Tian.
It is known that properness is independent of the K\"{a}hler class $\omega$ (see \cite{psc}). In \cite{psc} Tian has proved (under the assumption $\eta(X)=\{0\}$ \footnote{$\eta (X)$ denotes the Lie algebra of holomorphic vector fields.}) that \emph{$\nu_{\omega}$ is proper iff there is a K\"ahler Einstein metric in the class $\omega$}.
 
 In order to detect properness (conjecturally) one restricts attention to the subspace of \emph{Bergman metrics} inside $P(X,\omega)$ as these are \emph{dense} in $P(X,\omega)$ (see  \cite{tianberg}, \cite{ruan}, \cite{zel}, \cite{cat} ). The Bergman metrics are induced by the Kodaira embeddings furnished by large multiples of the polarization. More precisely, let $X\hookrightarrow \mathbb{C}P^{N}$ be the Kodaira embedding defined by a basis $\{S_0,S_1,\dots ,S_{N}\}$ of $H^{0}(X,L^{\otimes m})$.
There is a map $i_{\{S_{j}\}}:G \ra P(X,\omega)$ given as follows
\begin{align*}
 i_{\{S_{j}\}}(\sigma)\equiv \frac{1}{m}\varphi_{\sigma}=\frac{1}{m}\log \left({\sum_{k=0}^{N}||\sum_{l=0}^{N}\sigma_{k,l}S_l||^2}\right) \ .
\end{align*}
$\omega +\tp \dl\dlb\frac{1}{m}\varphi_{\sigma}$ are the Bergman metrics.
If $E_{\omega}$ denotes one of the energies we let $E_{\omega}(\sigma):= E_{\omega}(\frac{1}{m}\varphi_{\sigma})$.
The main issue is to let $\sigma=\lam_{t}$, an algebraic one parameter subgroup, and analyze the small $t$ asymptotics $\lim_{t\ra 0} E_{\omega}(\lam_{t})$. 
It is when we restrict the energies to the subspaces defined by the Bergman metrics 
that we make contact with Geometric Invariant Theory. 
 The first result in this direction is due to Tian who established the following theorem by exhibiting the K-energy map as the logarithm of a singular metric on a power of  the ample divisor on the moduli space of hypersurfaces of degree $d$ in $\cpn$. This is unquestionably the paradigm for all subsequent results in this area of investigation.
\ \\

 \noindent \textbf{Theorem} (G. Tian \cite{kenhyp})
 \emph{Let $Z_{f}$ be a normal degree $d\geq 2$ hypersurface in $\cpn$ , then $Z_{f}$ is stable if the K-energy is proper, and $Z_{f}$ is semistable if the K-energy is bounded from below.}
\ \\ 
\ \\
There is a related result in higher codimension, which was established independently by S. Zhang and the author.  Simon Donaldson has found an outstanding application of this theorem (see \cite{skdproj}).
\ \\
\ \\
\noindent \textbf{Theorem} (Zhang \cite{zhang}, Paul \cite{gacms})
\emph{Let $X$ be an n-dimensional subvariety of \cpn, and let $Z_{X}:= \{R_{X}=0 \}$ denote the associated hypersurface ($R_{X}$ is the $X$-resultant) . Let $\lambda(t)$ be an algebraic one parameter subgroup of  $G$. The weight of the action of $\lambda$ on $R_{X}$ is denoted by $F_{0}(\lambda)$. Then the following asymptotic expansion holds as $|t| \ra 0$}
 \begin{align*}
 d(n+1)F^0_{\omega}( \varphi_{\lambda (t)}) = F_{0}(\lambda)\log(|t|^2) +O(1).
\end{align*}

 In 1992 Ding and Tian (see \cite{dingtian} ) have studied the K-energy asymptotics along the \emph{algebraic potentials} $\vp_{\lambda(t)}$ when the limit is an almost Fano variety. They defined the \emph{generalized} Futaki invariant of a degeneration and proved that the \emph{sign} of this invariant (which is a rational number) is an obstruction to the existence of a K\"ahler Einstein metric in the class $-K_{X}$.  In 2002 (see \cite{skdtoric}) Simon Donaldson connected this invariant in an exciting way to Geometric Invariant Theory on the Hilbert scheme. Recently the author and G.Tian have proved the following theorem.
  \ \\
  \ \\
  \noindent \textbf{Theorem} (Paul, Tian \cite{cms2})
 \emph{Assume that $(X,L)$ moves in a  good family \footnote{See \cite{cms2}.} $\mathfrak{X}$. Then there is a function
$\Psi_{\mathfrak{X}}:G \rightarrow \mathbb{R}$  such that $-\infty \leq \Psi_{\mathfrak{X} }\leq C$  and an asymptotic expansion}
\begin{align*}
d(n+1){ \nu_{\omega}(\varphi_{\lambda(t)})-\Psi_{\mathfrak{X}}({\lambda(t)})= 
 F_{1}(\lambda)\log(|t|^2)+O(1)} \ \mbox{as}\ |t| \rightarrow 0. \quad 
\end{align*}
\emph{Moreover  $\Psi_{\mathfrak{X}}(\lambda(t))= \psi(\lambda)\log(|t|^{2})+O(1)$  where  $\psi({\lambda}) \in \mathbb{Q}_{\geq 0} $ and  $\psi({\lambda}) \in \mathbb{Q}_{+} $ if and only if $X^{\lambda(0)}$ (the limit cycle  \footnote {See \cite{sopv} pg. 61.} of $X$ under $\lambda$ ) has a component of multiplicity greater than one.  $F_{1}(\lambda)$ is the generalized Futaki invariant  \footnote{See \cite{skdtoric} and \cite{dingtian}.}of the degeneration $\lambda$, and $O(1)$ denotes any quantity which is bounded as $|t|\rightarrow 0$.}
\ \\
\ \\
 In their study of the K\"ahler Ricci flow on K\"ahler Einstein manifolds X.X. Chen and G.Tian introduced a set of new energy functionals  $E_{\om,l}$ ($l=0,1,2,\dots n$) which monotonically decrease along the flow (under a positivity hypothesis). These have received much attention in the recent literature.    
   It is now known that these energies are bounded from below in the presence of a K\"ahler Einstein metric. Therefore it is reasonable to attempt to connect these energies to Geometric Invariant Theory
 in the spirit of the above theorems. Unfortunately this does not seem possible, however in this paper we are able to modify the $E_{\om,l}$  so that the asymptotic results go through. We call these new energies simply \emph{higher energies}.
\begin{center}{\sc{Acknowledgements}}\end{center}
The author would like to thank Xiu Xiong Chen for many helpful conversations.
 
  
  \begin{center}{\S 1 \sc{Statement of Results}}\end{center}
 Let $X\hookrightarrow \cpn$ be a projective variety. We are concerned with a certain procedure the rough expression of which is the following.
\begin{flushleft}
\emph{Step one}. Given $\lambda:\mathbb{C}^*\ra G$ take the coefficient $F_{l}(\lambda)$ of $m^{-l}$ in the expansion of the normalized weight of the action of $\lambda$ on the $mth$ Hilbert point of $X$.\\
\ \\
\emph{Step two}. Choose a virtual bundle $\xi$ such that the determinant of its direct image has weight $F_{l}(\lambda)$ with respect to the action of $\lambda$.\\
 \ \\
 \emph{Step three}. Define an energy functional $F_{l,\om}$ as the transgression of the Riemann Roch Hirzebruch integrand with respect to $\xi$. General theory then exhibits this energy as a singular metric on the sheaf obtained in step two.\\
 \ \\
\emph{Conclusion}. Algebro geometric energy asymptotics of the functional $F_{l,\om}$.
 \end{flushleft}

More precisely let $\mathfrak{X}\overset{f}\ra S$ be a flat family of subschemes of $\cpn$. Let  $\mathcal{L}\in \mbox{Pic}(\mathfrak{X})$ be relatively very ample. $\mathbb{Q}[\mathcal{L}]$ denotes the subring of $K^{\circ}(\mathfrak{X})\otimes_{\mathbb{Z}}\mathbb{Q}$ generated by $\mathcal{L}$.  Let $\Delta$ denote the discriminant locus of $f$. That is, $X_{z}$ is $C^{\infty}$ for all $z \in S\setminus \Delta$.
 Our purpose is to study two maps
\begin{align*}
& i) \ \tau_{X}:\mathbb{Q}[\mathcal{L}]\times P_{\om}(X_{z})\rightarrow \mathbb{R} \quad  z\in S\setminus \Delta \\
&\tau_{X}(\xi,\vp):=\left(\frac{\sqrt{-1}}{2\pi}\right)^{n}\int_{X}\int_{0}^{1}\frac{\dl}{\dl b}\left(\mbox{Td}(R_{g_{t}}+b\frac{\dl g_{t}}{\dl t}g_{t}^{-1})\mbox{Ch}(F^{\xi}_{H_{t}}+b\frac{\dl H_{t}}{\dl t}H_{t}^{-1})\right)_{|_{b=0}}dt \\
 &\mbox{Where}\\
 &g_{t}=\omega+\frac{\sqrt{-1}}{2\pi}\dl\dlb \varphi_{t}, \ \mbox{a family of K\"ahler metrics on} \ X,\\
&H_{t} =e^{-\varphi_{t}}h, \ h \ \mbox{a Hermitian metric on}\  \mathcal{L},\\
&F^{\xi}_{H_{t}}, \ \ \mbox{The curvature of $\xi$ with respect to $H_{t}$},\\
&\omega=\frac{-\sqrt{-1}}{2\pi}\dl\dlb \mbox{log}(h),\ \mbox{the curvature of} \ h,\\
&R_{g}:= \dlb \{\dl(g) g^{-1}\}\in C^{\infty}(\Lambda^{(1,1)}\otimes \mbox{End}(T^{1,0}_{X})),\ \mbox{the full K\"ahler curvature tensor.}\\
& \mbox{and}\\
 & ii)\  \mathbf{Det}:\mathbb{Q}[\mathcal{L}]\rightarrow \mbox{Pic}(S), \
   \mathbf{Det}(\xi):= \mbox{det}(Rf_{*}(\xi))\ .
 \end{align*}
 For the virtual bundle in step two we make the following choice:
  \begin{align}
\mathcal{E}_{l}(m):=\sum_{\{0\leq j \leq l\}}\sum_{\{0\leq i \leq n+1-j\}}q_{n+1-j}(m)(-1)^{i}\binom{n+1-j}{i}\mathcal{L}^{m+i} \in \mathbb{Q}[\mathcal{L}].
\end{align}
 where  $0\leq l \leq n+1$. The coefficients $q_{n+1-j}(m)\in \mathbb{Q}$ are defined below. 
 
 \begin{theorem}
 
 There are invertible sheaves $\mathcal{L}_{l}$ \emph{ ($0\leq l \leq n+1$)} on $ \mathfrak{Hilb}_{\mathbb{P}^{N}}^{\chi}(\mathbb{C})$ such that for any family  $\mathfrak{X}\overset{f}\ra S$ as above we have
 \begin{align*}
& i) \ \mathbf{Det}(\mathcal{E}_{l}(m))\cong g^{*}(\mathcal{L}_{l}). \
\mbox{Consequently  $\mathbf{Det}(\mathcal{E}_{l}(m))$ is independent of $m$}.\\
\ \\
&\mbox{Assume that  $S$ is proper, and that $G$ acts on the data,}\\
&\mbox{then we can state the following.}\\
\ \\
& ii)\ \mbox{For any one parameter subgroup  $\lambda$  of $G$ and $z\in  \mathfrak{Hilb}_{\mathbb{P}^{N}}^{P}(\mathbb{C})$}\\
&\mbox{let $w_{\lambda}(z)$ denote the weight  of the $\mathbb{C}^{*}$ action on 
${{\mathcal{L}^{\vee}}_{l}}|_{z_{0}}$  where  $z_{0}=\lambda(0)z$,}\\
&\mbox{then} \ {w_{\lambda}(z)=F_{l}(\lambda)}\ . \\
\ \\
& \mbox{Assume that GRR is valid for the map $\mathfrak{X}\overset{f}\ra S$. Then}\\
\ \\
&iii) \ c_{1}(g^{*}\mathcal{L}_{l})= \sum_{n+1-l\leq k \leq n+1}f_{*}\{ \frac{c_{l,k}}{k!}\mbox{Td}_{n+1-k}(f)\om^{k}\}\ .
 \end{align*}
\end{theorem}
Next we introduce our new energy functionals. These simultaneously generalize the Aubin energy, Mabuchi's K-energy, and (modulo lower order terms) the Chen Tian energy functionals. General theory shows that they are all path independent.
  \begin{definition} (Higher Energy Functionals)
 \begin{align}
  F_{l,\om}(\vp):= \tau_{X}(\mathcal{E}_{l}(m),\vp)
\end{align}
\end{definition}
 In the statement of the next result we assume that $X=X_{z}$ is a generic member in a smooth family $\mathfrak{X}\overset{f}\ra S$.
 \begin{theorem}
 Let $\lambda$ be a smooth degeneration of $X$.  Then for all $0\leq l \leq n+1$ there is an asymptotic expansion
\begin{align}
  F_{l,\om}(\varphi_{\lambda(t)})= 
 F_{l}(\lambda)\log(|t|^2)+O(1) \ \mbox{as}\ |t| \rightarrow 0. 
\end{align}
Where  $O(1)$ denotes any quantity which is bounded as $|t|\rightarrow 0$. Moreover, we have
 \begin{align*}
& F_{0,\om}(\varphi)= F^{0}_{\om}(\vp) \ (\mbox{Aubin's Energy})\\
&\ \\
 &F_{1,\om}(\vp)=\nu_{\om}(\vp) \ (\mbox{Mabuchi's K-energy})\ .
 \end{align*}
 \end{theorem}
\begin{corollary}
If $F_{l,\om}(\varphi)$ is bounded from below then for all smooth degenerations $\lambda$ we have that $F_{l}(\lambda)\leq 0$.
\end{corollary}
 Recall that  $E_{\om,l}$ is given by the following expression (see \cite{chentianricci} for more details).
 \begin{align*}
\frac{ d E_{\om,l}}{dt}:=&\frac{l+1}{V}\int_{X}\Delta_{\vp}\left(\frac{\dl \vp}{\dl t}\right)\mbox{Ric}(\om_{\vp})^{l}\wedge \om^{n-l}_{\vp}-
\frac{n-l}{V}\int_{X}\frac{\dl \vp}{\dl t}(\mbox{Ric}(\om_{\vp})^{l+1}-\om^{l+1}_{\vp})\wedge \om^{n-l-1}_{\vp}.
 \end{align*}
In particular when $l=1$ we have that
  \begin{align*}
\frac{dE_{\om,1}}{dt}=\frac{2n}{V}\int_{M}\Delta_{\vpt}\dot\vpt\mbox{Ric}(\vpt)\wedge\om_{\vpt}^{n-1}-\frac{n(n-1)}{V}\int_{M}\dot\vpt\left(\mbox{Ric}({\vpt})^{2}-\om_{\vpt}^{2}\right)\wedge\om_{\vpt}^{n-2}.
\end{align*}
 In this paper we limit our study to the first new energy functional $F_{2,\om}$. We will take up the further study of this functional in a sequel to this paper.
  Next we state the following proposition which compares $E_{1,\om}$ with $F_{2,\om}$.  
 \begin{proposition}
 \begin{align*}
\frac{dF_{2,\om}(\vpt)}{dt} =&\frac{2n}{V}\int_{M}\Delta_{\varphi_{t}}\dot\varphi_{t}\mbox{\emph{Ric}}(\varphi_{t})\wedge\om_{\vpt}^{n-1}
-\frac{3n(n-1)}{2V}\int_{M}\dot\vpt\left(\mbox{\emph{Ric}}(\vpt)^{2}-\alpha(g)\om^{2}_{\vpt}\right)\wedge\om_{\vpt}^{n-2}\\
&-\frac{1}{V}\int_{M}\dot\vpt\left( ||\mbox{R}(\vpt)||^{2} - ||\mbox{\emph{Ric}}(\vpt)||^{2}- \beta(g) \right)\om_{\vpt}^{n}-3{\mu}\frac{d\nu_{\om}}{dt}
\end{align*}
\end{proposition}
\noindent $||R||^{2}$ and $||\mbox{Ric}||^{2}$ denote the square norm of the full curvature tensor and Ricci curvature respectively. 
\begin{align*}
&\alpha(g):=\frac{1}{V}\int_{M}\mbox{Ric}_{g}^{2}\wedge\frac{\om_{g}^{n-2}}{n!}, \ \beta(g):=\frac{1}{V}\int_{M}\left(||\mbox{R}||^{2}- ||\mbox{Ric}_{g}||^{2}\right)\frac{\om^{n}}{n!}  \\
\end{align*}

 
  \begin{center}{\S 1 \sc{Hilbert Points and Chow Forms}}\end{center}
 
Let $(X,\mathcal{L})$ be a polarized algebraic variety. Assume that $\mathcal{L}$ is very ample with associated embedding
 \begin{align*}
 X\underset{\varphi_{\mathcal{L}}}\longrightarrow \mathbb{P}(H^{0}(X,\mathcal{L})^{*}).
 \end{align*}
 Fix  an isomorphism
 \begin{align*}
 \sigma:H^{0}(X,\mathcal{L})^{*}\underset{\cong}\longrightarrow \mathbb{C}^{N+1}\ .
 \end{align*}
 In this way we consider $X$ embedded in $\cpn$.
 Let $m\in \mathbb{Z}$ be a large positive integer.
 Then there is a surjection\footnote{$\mathbf{S}^{m}$ denotes the $m$th symmetric power operator.}
 \begin{align*}
 \psi_{X,m}:\mathbf{S}^{m}({\mathbb{C}^{N+1}})^{*}\longrightarrow H^{0}(X,\mathcal{O}(m))\rightarrow 0\ .
 \end{align*}
 Let $\chi(m)=\chi(X,\mathcal{O}(m))=h^{0}(X,\mathcal{O}(m))$. It is a deep fact (see \cite{eckart}) that there is an integer $m(P)$ depending \emph{only} on the Hilbert polynomial $P$ such that for all $m\geq m(P)$, the kernel of $ \psi_{X,m}$ 
 \begin{align*}
 \mbox{Ker}(\psi_{X,m})\in \mbox{G}(\chi(m),\mathbf{S}^{m}({\mathbb{C}^{N+1}})^{*}) 
 \end{align*}\footnote{The Grassmannian of $\chi(m)$ dimensional \emph{quotients} of $\mathbf{S}^{m}({\mathbb{C}^{N+1}})^{*}$.}
 completely determines $X$. In other words,  the \emph{entire} homogeneous (saturated) ideal can be recovered from its $m$th graded piece.
 We have the Pl\"ucker embedding
\begin{align*}
 \mathcal{P}:\mbox{G}(\chi(m),\mathbf{S}^{m}({\mathbb{C}^{N+1}})^{*}) \rightarrow \mathbb{P}\left(\bigwedge^{d_{m}-\chi(m)}\mathbf{S}^{m}(\mathbb{C}^{N+1})^{*}\right).
 \end{align*}
Next we consider the canonical nonsingular  pairing 
 \begin{align*}
 \bigwedge^{d_{m}-\chi(m)}\mathbf{S}^{m}(\mathbb{C}^{N+1})^{*}\otimes \bigwedge ^{\chi(m)}\mathbf{S}^{m}(\mathbb{C}^{N+1})^{*}\longrightarrow \mbox{\bf{det}}(\mathbf{S}^{m}(\mathbb{C}^{N+1})^{*})\ .
 \end{align*}
 This induces a natural isomorphism
 \begin{align*}
 \mathbb{P}\left(\bigwedge^{d_{m}-\chi(m)}\mathbf{S}^{m}(\mathbb{C}^{N+1})^{*}\right)\overset{\iota}\cong \mathbb{P}\left(\bigwedge^{\chi(m)}\mathbf{S}^{m}(\mathbb{C}^{N+1})\right).
 \end{align*}
 Combining this identification with the Pl\"ucker embedding we associate to $\mbox{Ker}(\psi_{X,m})$ a unique point, called (following Gieseker) the $m$th Hilbert Point
 \begin{align*}
\mbox{Hilb}_{m}(X):= \iota(\mathcal{P}(\mbox{Ker}(\psi_{X,m})))\in  \mathbb{P}\left(\bigwedge^{\chi(m)}\mathbf{S}^{m}(\mathbb{C}^{N+1})\right).
\end{align*}
Let $\mathbf{I}=(i_{0},i_{1},\dots,i_{N})$ be a multiindex with $|\mathbf{I}|:=i_{0}+i_{1}+\dots +i_{N}=m$, $i_{j}\in \mathbb{N}$. Let $e_{0},e_{1},\dots, e_{N}$ be the standard basis of $\mathbb{C}^{N+1}$, and $z_{0},z_{1},\dots,z_{N}$ be the dual basis of linear forms. Consider the monomials
 $M_{\mathbf{I}}:= e_{0}^{i_{0}}e_{1}^{i_{1}}\dots e_{N}^{i_{N}}$ and $M^{*}_{\mathbf{I}}:= z_{0}^{i_{0}}z_{1}^{i_{1}}\dots z_{N}^{i_{N}}$. Fix a basis $\{f_{1},\dots, f_{\chi(m)}\}$ of $H^{0}(X,\mathcal{O}(m))$. 
 Then 
\begin{align*}
\begin{split}
\wedge^{\chi(m)}\psi_{X,m}(M^{*}_{{\mathbf{I}}_{j_{1}}}\wedge \dots \wedge M^{*}_{{\mathbf{I}}_{j_{\chi(m)}}})=&\psi_{X,m}(j_{1},\dots , j_{\chi(m)})f_{1}\wedge \dots \wedge  f_{\chi(m)} \\
 &\psi_{X,m}(j_{1},\dots , j_{\chi(m)})\in \mathbb{C}\ .
\end{split}
\end{align*}
Then in homogeneous coordinates we can write
\begin{align*}
\wedge^{\chi(m)}\psi_{X,m}=\sum_{\{{\mathbf{I}}_{j_{1}},\dots,{\mathbf{I}}_{j_{\chi(m)}}\}}\psi_{X,m}(j_{1},\dots , j_{\chi(m)})    M_{{\mathbf{I}}_{j_{1}}}\wedge \dots \wedge M_{{\mathbf{I}}_{j_{\chi(m)}}}.
 \end{align*}
  Let $\lambda:\mathbb{C}^{*}\rightarrow G$ be an algebraic one parameter subgroup. We may assume that $\lambda$ has been diagonalized on the standard basis $\{e_{0},e_{1},\dots, e_{N}\}$. Explicitly, we assume that there are  $r_{i}\in \mathbb{Z}$ such that
 \begin{align*}
 \lambda(t)e_{j}=t^{r_{j}}e_{j} .
 \end{align*}
Define the weight of $\lambda$ on the monomial $M_{\mathbf{I}}$ by
\begin{align*}
w_{\lambda}(M_{\mathbf{I}}):= r_{0}i_{0}+r_{1}i_{1}+\dots +r_{N}i_{N}.
\end{align*}
 
\begin{definition} (Gieseker \cite{globmod})
The  \emph{\textbf{weight}} of the $m$th Hilbert point of $X$ is the integer 
\begin{align*}
w_{\lambda}(m):=\overset{\mbox{\emph{Min}}}{\{{\mathbf{I}}_{j_{1}},\dots,{\mathbf{I}}_{j_{\chi(m)}}\}} \left(\sum_{1\leq k\leq \chi(m)}w_{\lambda}(M_{{\mathbf{I}}_{j_{k}}})|\psi_{X,m}(j_{1},\dots , j_{\chi(m)})\neq 0\right)\ .
\end{align*}
\end{definition}
  Let $X \subset \cpn$ be an $n$ dimensional irreducible subvariety of $\cpn$ with degree $d$, then the Chow form,
or associated hypersurface to $X$ is defined by
\begin{align*}\label{chow}
Z_{X}:= \{L \in \mathbb{G}:= \mathbb{G}(N-n-1,\mathbb{C}P^{N}): L\cap X \neq \emptyset\}.
\end{align*}
It is easy to see that $Z_{X}$ is an \emph{irreducible} hypersurface (of degree $d$) in $\mathbb{G}$.
Since the homogeneous coordinate ring of the grassmannian is a UFD, any codimension one subvariety with
degree $d$ is given by the vanishing of a section  $R_{X}$ of the homogeneous coordinate ring\footnote{See \cite{tableaux}  pg. 140 exercise 7.}
 \begin{align*}
 \{\ R_{X}=0\ \}= Z_{X}  \ ; \  R_{X} \in \mathbf{P}H^{0}(\mathbb{G},\mathcal{O}(d)).
 \end{align*}
$R_{X}$ is confounded with $Z_{X}$.
Following \cite{ksz} we can be more concrete as follows.
Let $M_{n+1, N+1}^{0}(\mathbb{C})$ be the (Zariski open and dense) submanifold
of the vector space of $M_{n+1, N+1}(\mathbb{C}) $ matrices of full rank.
We have the canonical projection
\begin{align*}
p:M_{n+1, N+1}^{0}(\mathbb{C})\rightarrow \gr ,
\end{align*}
defined  by taking the kernel of the linear transformation. This map is dominant, so the closure of the  preimage
\begin{align*}
\overline{p^{-1}(Z_{X})}\subset \overline{M_{n+1, N+1}^{0}(\mathbb{C})}=
M_{n+1, N+1}(\mathbb{C})
\end{align*}
is also an irreducible hypersurface of degree $d$ in $M_{n+1, N+1}(\mathbb{C})$.
Therefore, there is a unique\footnote{Unique up to scaling.} (symmetric multihomogeneous) polynomial (which will also be denoted by $R_{X}$)
such that 
\begin{align*}
Z:= \overline{p^{-1}(Z_{X})}= \{R_{X}(w_{ij})=0\} \ ;  \ R_{X}(w_{ij})\in \mathcal{P}^{d}[M_{n+1, N+1}(\mathbb{C})].
 \end{align*}
   There is  (in principal) an \emph{explicit formula} for the polynomial $R_{X}(w_{ij})$, which is essentially due to Cayley in his remarkable 1848 note \cite{cay} on resultants.  The modern formulation of these ideas are due to Grothendieck, Knudsen, and Mumford. See \cite{detdiv}. 
{\textbf{Theorem}} \emph{(Cayley, Grothendieck, Knudsen, Mumford)}\\
\emph{There is a canonical isomorphism of one dimensional vector spaces}
\begin{align*}
&\Delta^{n+1}{\bf{det}}(H^{0}(X,\mathcal{O}(m)))={\bigotimes_{i=0}^{n+1}{\bf{det}}(H^{0}(X,\mathcal{O}(m+i)))^{(-1)^{i+1}\binom{n+1}{i}}\cong \mathbb{C}R_{X}^{(-1)^{n+1}}}.\\
&\Delta \ \mbox{denotes the first forward difference operator.}
\end{align*}
  
  It follows from this that the weight is a  $w_{\lambda}(m)$ is a \emph{polynomial} in $m$
 \begin{align*}
w_{\lambda}(m)=a_{n+1}(\lambda)m^{n+1}+a_{n}(\lambda)m^{n}+O(m^{n-1}).
\end{align*}
Let $\chi(m)$ be the Hilbert polynomial of $X$. Following Gieseker, we consider the ratio
\begin{align*}
\frac{w_{\lambda}(m)}{m\chi(m)} \ .
\end{align*}
As in Donaldson \cite{skdtoric} we consider the coefficients of $m^{-l}$ in the expansion
\begin{align*}
&\frac{w_{\lambda}(m)}{m\chi(m)} =F_{0}(\lambda)+F_{1}(\lambda)\frac{1}{m}+\dots +F_{l}(\lambda)\frac{1}{m^{l}}+\dots\\
\ \\
&F_{l}(\lambda)= c_{l,n+1}a_{n+1}(\lambda)+c_{l,n}a_{n}(\lambda)+c_{l,n-1}a_{n-1}(\lambda)+\dots + c_{l,n+1-l}a_{n+1-l}(\lambda).
\end{align*}
The $c_{l,j}$ are all rational functions of the coefficients of the Hilbert polynomial $\chi$. \\

The relationship between Hilbert Points and Chow forms extends to the relative situation as well. Let $f:\mathfrak{X}\ra S$ be a flat morphism of projective varieties. 
\begin{diagram}
\mathfrak{X}&\rInto^{\iota}&\mathbb{P}(f_{*}\mathcal{L}^{\vee})\cong S\times \mathbb{P}^{N}\\
\dTo^{f}&\ldTo_{\pi}&\\
  S&
\end{diagram}
We assume that $\mathcal{L}$ is a relatively ample line bundle on $\mathfrak{X}$ with respect to the map $f$.
The isomorphism $\mathbb{P}(f_{*}\mathcal{L}^{\vee})\cong S\times \mathbb{P}^{N}$ is equivalent to the existence of an invertible sheaf $\mathcal{A}$ on $S$ such that 
\begin{align*}
 f_{*}\mathcal{L}\cong \underset{N+1}{\underbrace{\bigoplus}} \mathcal{A}\ .  
\end{align*}
Let $\chi$ denote the Hilbert polynomial of the fibers. Let 
$ \mathfrak{Hilb}_{\mathbb{P}^{N}}^{\chi}(\mathbb{C})$ denote the \emph{Hilbert scheme}. For $m$ large enough there is a map $\vp_{m}$ from $S$ into  
$ \mathfrak{Hilb}_{\mathbb{P}^{N}}^{\chi}(\mathbb{C})$
\begin{align*}
\begin{split}
&\varphi_{m}:S\rightarrow \mathfrak{Hilb}_{\mathbb{P}^{N}}^{\chi}(\mathbb{C})\hookrightarrow \mathbb{P}^{N(m)}:= \mathbb{P}(\bigwedge^{\chi(m)}\mathbf{S}^{m}
(\underset{N+1}{\underbrace{\bigoplus}} \mathbb{C})) .\\ 
&\mbox{ Pulling back $\mathcal{O}_{\mathbb{P}^{N(m)}}(1)$ to $S$ via $\vp_{m}$ gives the isomorphism}\\
&\varphi_{m}^{*}\mathcal{O}_{\mathbb{P}^{N(m)}}(1)\cong \mbox{det}(f_{*}\mathcal{L}^{m})\otimes \mbox{det}(f_{*}\mathcal{L})^{\frac{-m\chi(m)}{N+1}} \ .
\end{split}
\end{align*}
Therefore the  appropriate generalization of Hilbert points to families is the following invertible sheaf on $S$
\begin{align*}
{\mbox{Hilb}_{m}(\mathfrak{X}\setminus S) :=  \mbox{det}(f_{*}\mathcal{L}^{m})\otimes \mbox{det}(f_{*}\mathcal{L})^{\frac{-m\chi(m)}{N+1}}}\ .
\end{align*}
 
Let $\mathfrak{C}(n,d; \cpn)$ denote the {Chow Variety}  of dimension $n$  and degree $d$ algebraic cycles inside $\cpn$.    
 There is a morphism  $\Delta^{n+1}$ from the {Hilbert scheme}  to the Chow variety  (see \cite{git} and \cite{fogarty}) which sends a subscheme $\mathcal{I} $ of $\cpn$ with Hilbert polynomial $\chi$ to the Chow form of the top dimensional component of its underlying cycle.    
 We now have the sequence of maps
 \begin{align*}
S\overset{\varphi_{m}}\longrightarrow \mathfrak{Hilb}_{\mathbb{P}^{N}}^{\chi}(\mathbb{C})\overset{\Delta^{n+1}}\longrightarrow \mathfrak{C}(n,d; \cpn) \overset{\iota}\hookrightarrow 
\mathbb{P}(H^{0}(\mathbb{G},\mathcal{O}(d))).
\end{align*}
Let $\mathcal{O}(1)$ denote the hyperplane line on $\mathbb{P}(H^{0}(\mathbb{G},\mathcal{O}(d)))$.
Then we define the \emph{Chow form of the map} $\mathfrak{X}\overset{f}\rightarrow S$ to be the invertible sheaf on $S$
\begin{align*}
\mbox{Chow}(\mathfrak{X}\overset{f}\rightarrow S ):= \varphi_{m}^{*}{\Delta^{n+1}}^{*}\iota^{*}\mathcal{O}(1).
\end{align*}
The extension to families of the relationship between Chow forms and Hilbert points can be stated as follows \\
\ \\\noindent  \emph{On the base $S$ there is a canonical isomorphism of invertible sheaves} 

 \begin{align*}
 \mbox{Chow}(\mathfrak{X}\overset{f}\rightarrow S)\otimes \mbox{det}(f_{*}\mathcal{L})^{\frac{d(n+1)}{N+1}} \cong \Delta^{n+1}\mbox{det}(f_{*}\mathcal{L}^{m})\ .
 \end{align*}
This immediately implies that the following expansion holds (see \cite{detdiv}, in particular Theorem 4, for a complete discussion).\newpage
\noindent \textbf{Theorem} (\emph{Knudsen, Mumford})\\
\noindent \emph{There are invertible sheaves $\mathcal{M}_{0},\mathcal{M}_{1},\dots, \mathcal{M}_{n+1} $ on $S$ and a canonical and functorial isomorphism}:
 \begin{align*}
 &\mbox{det}(f_{*}\mathcal{L}^{m})\cong \bigotimes _{j=0}^{n+1}\mathcal{M}_{j}^{\binom{m}{j}}\ .
  \end{align*}
 
 \begin{center}{\S 2 \sc{Higher Sheaves}}\end{center}
Let $\omega_{\mathfrak{X}/ S}:= K_{\mathfrak{X}}\otimes f^{*}({K_{S}}^{-1})$ denote  the relative canonical bundle. In \cite{kenhyp} Tian introduced the following invertible sheaf on $S$ , which one calls the \emph{CM polarization}.
\begin{align*}
\mbox{L}_ {CM}:= \mbox{det}\left(f_{*}2^{-(n+1)}\left((n+1)(\omega_{\mathfrak{X}/ S}^{-1}- \omega_{\mathfrak{X}/ S})(\mathcal{L}-\mathcal{L}^{-1})^{n}-\mu(\mathcal{L}-\mathcal{L}^{-1})^{n+1}\right)\right)^{-1}.
\end{align*}
We should emphasize that the CM polarization inspired though it is by the G.I.T. approach to moduli, naturally appeared in the setting of analysis and differential geometry.   For many reasons  it became necessary to have a \emph{purely algebraic} construction of this sheaf. The following result extends the definition of this sheaf to \emph{any} flat family of subschemes of some $\cpn$. In particular to the universal family over the Hilbert Scheme.  This extension is based on the theory of determinants of Cayley, Grothendieck,  Knudsen and Mumford. \\
\noindent \textbf{Theorem} (Paul, Tian \cite{cms1}) \ \\
\emph{ We make the following assumptions\\
\ \\
i)  $\mathfrak{X}\overset{f}\rightarrow S$ is a flat, proper, local complete intersection morphism\footnote{See \cite{fulton}.}.\\
\ \\
ii) The map $\mbox{Pic}(S)\overset{c_{1}}\rightarrow \mbox{H}^{\ 2}(S,\mathbb{Z})$ is injective.\\
\ \\
 iii) $\mathcal{L}$ is relatively ample on $\mathfrak{X}$ and  $\mathbb{P}(f_{*}\mathcal{L})\cong S\times \mathbb{P}^{N}$.\\
  \ \\
  Then there is a canonical and functorial isomorphism of sheaves on $S$ }
\begin{align*}
{\mbox{L}_{CM}\cong \{\mbox{Chow}(\mathfrak{X})\otimes \mbox{det}(f_{*}\mathcal{L})^{\frac{d(n+1)}{N+1}}\}^{n(n+1)+\mu}\otimes\mathcal{M}_{n}^{-2(n+1)}} \ .
\end{align*}
As the reader will see, on the right hand side of this isomorphism is the sheaf $\mathcal{L}_{1}$.
Let $f$ be any numerical function,  recall that the forward difference of $f$ is defined as follows.
\begin{align*}
&\Delta f (m):=f(m+1)-f(m)
\end{align*}
{Inductively we set}
\begin{align*}
\Delta^{k+1}f(m):= \Delta^{k}f(m+1)-\Delta^{k}f(m).\\
\end{align*}
 Let $f_{l}(m):=m^{l}$. Then we define polynomials $P_{k,l}(m)$
\begin{align*}
P_{k,l}(m):=\Delta^{k}f_{l}(m) \ . 
\end{align*}
It is easy to see that
\begin{align*}
P_{k,l}(m)=\sum_{0\leq j \leq k}(-1)^{j+1}\binom{k}{j}(m+j)^{l}\ .
\end{align*}
It is not difficult to verify that
\begin{align*}
P_{k,l}(m)=
\begin{cases}
(-1)^{k}k! , & \mbox{if} \ k=l \\
0,& \mbox{if} \  l<k \ .
\end{cases}
\end{align*}
In general, $P_{k,k+d}(m)$ is a polynomial in $m$ of degree $d$.
Given $0\leq l \leq n+1$ let 
\begin{align*}
(q_{n+1}(m), q_{n}(m), q_{n-1}(m), \dots, q_{n+1-l}(m))
\end{align*}
be the unique solution to the equation
\begin{align*}
\begin{pmatrix} P_{n+1,n+1}(m)&P_{n,n+1}(m)&\dots & \dots & P_{n+1-l,n+1}(m)\\
      0&P_{n,n}(m)& P_{n-1,n}(m)&\dots & P_{n+1-l,n}(m)  \\
       0& 0&P_{n-1,n-1}(m)  &\dots& P_{n+1-l,n-1}(m)\\
       0&  0& 0& \dots &\dots  \\ 
       \dots&\dots&\dots&\dots\\
         0&  0& 0& \dots&  P_{n+1-l,n+1-l}(m)
       \end{pmatrix}
\begin{pmatrix} q_{n+1}(m)\\
q_{n}(m)\\
q_{n-1}(m)\\
\dots \\
\dots\\
q_{n+1-l}(m)
\end{pmatrix}
=
 \begin{pmatrix} c_{l,n+1}\\
c_{l,n}\\
c_{l,n-1}\\
\dots \\
\dots\\
c_{l,n+1-l}
\end{pmatrix}
\end{align*}
 
\begin{definition} (Higher Sheaves)
Let $\mathfrak{X}\overset{f}\rightarrow S$ be a flat family of subschemes of $\cpn$. Let $\mathcal{L}$ be ample with respect to $f$. We define sheaves $\mathcal{L}_{l}$ for all $l=0,1,2,\dots,n+1$ on $S$ as follows
 \begin{align}
\mathcal{L}_{l}:=\bigotimes_{k=n+1-l}^{n+1}{\mathcal{M}_{k}}^{\frac{1}{k!}\sum_{0\leq j\leq k-1}(-1)^{j+1}\sigma_{j}(1,2,\dots,k-1)c_{l.k-j}}
\end{align}
Where the $\mathcal{M}_{k} ,\ 0\leq k \leq n+1$ are the coefficients in the Cayley, Grothendieck, Knudsen, Mumford expansion.
 \end{definition}
 Observe that we have the following canonical isomorphisms
\begin{align*}
&\mathcal{L}_{0}\cong \mbox{Chow}(\mathfrak{X}\overset{f}\rightarrow S)\otimes \mbox{det}(f_{*}\mathcal{L})^{\frac{d(n+1)}{N+1}} \\
&\mathcal{L}_{1}\cong \mbox{L}_{CM}\ .
 \end{align*}
Recall that the Hilbert point of the family $\mathfrak{X}\overset{f}\rightarrow S$ is the invertible sheaf
 \begin{align*}
&\mbox{Hilb}_{m}(\mathfrak{X}\setminus S) :=  \mbox{det}(f_{*}\mathcal{L}^{m})\otimes  \mathcal{A}^{\frac{-m\chi(m)}{N+1}}\\
&\mathcal{A}:= \mbox{det}(f_{*}\mathcal{L})\ .
\end{align*}
Then we have the following proposition.
\begin{proposition} For all $l=0,1,2 \dots, n+1$ we have
\begin{align}\label{sum0}
\bigotimes_{0\leq p\leq l}\bigotimes_{0\leq i \leq n+1-p} \mbox{\emph{Hilb}}_{m}(\mathfrak{X}\setminus S)^{(-1)^{i}q_{n+1-p}(m)\binom{n+1-p}{i}} \cong  \mathcal{L}_{l}
\end{align}
 
\end{proposition} 
\emph{Proof}\\
Writing out the left hand side of (\ref{sum0}) gives 
\begin{align*}
&\bigotimes_{0\leq p\leq l}\bigotimes_{0\leq i \leq n+1-p}\mathbf{det}(f_{*}\mathcal{L}^{m+i})^{(-1)^{i}q_{n+1-p}(m)\binom{n+1-p}{i}}\otimes \mathcal{A}^{\frac{1}{N+1}(-1)^{i+1}(m+i)\chi(m+i)q_{n+1-p}(m)\binom{n+1-p}{i}}\ . \\
\end{align*}
The exponent of $\mathcal{A}$  satisfies the following
\begin{align}\label{sum}
\sum_{0\leq p\leq l}\sum_{0\leq i\leq n+1-p}(-1)^{i+1}(m+i)\chi(m+i)q_{n+1-p}(m)\binom{n+1-p}{i}=\begin{cases}1,&l=0\\
0,& l>0\ .
\end{cases}
\end{align}
 
 To see this we first write
 \begin{align*}
 (m+i)\chi(m+i)=b_{n}(m+i)^{n+1}+b_{n-1}(m+i)^{n}+\dots + b_{j}(m+i)^{j+1}+\dots
\end{align*}
Then the left hand side of (\ref{sum}) is given by
\begin{align}\label{sum2}
\sum_{0\leq p\leq l}\sum_{0\leq i\leq n+1-p}\sum_{0\leq j \leq n}(-1)^{i+1}  q_{n+1-p}(m)b_{j}(m+i)^{j+1}\binom{n+1-p}{i}
\end{align}
Recall that we have defined the polynomials $P_{n+1-p,j+1}(m)$ by the formula
\begin{align}\label{P}
P_{n+1-p,\ j+1}(m)= \sum_{0\leq i \leq n+1-p}(-1)^{i+1}(m+i)^{j+1}\binom{n+1-p}{i}\ .
\end{align}
Substituting (\ref{P}) into (\ref{sum2}), switching the order of summation  and appealing to the definiton of the $q_{k}(m)$ gives
\begin{align}\label{sum3}
\sum_{n-l\leq j \leq n}\sum_{n-j\leq p\leq l}q_{n+1-p}(m)P_{n+1-p,\ j+1}(m)b_{j}=\sum_{n-l\leq j \leq n}b_{j}c_{l,\ j+1}.
\end{align}
By definition of the $c_{l,k}$ the right hand side of (\ref{sum3}) is the coefficient of ${m}^{-l}$ in the expansion of
\begin{align*}
\frac{m\chi(m)}{m\chi(m)}\equiv 1\ .
\end{align*}
From now on we will assume that $l>0$. With this assumption we have
\begin{align}
 &\bigotimes_{0\leq p\leq l}\bigotimes_{0\leq i \leq n+1-p} \mbox{{Hilb}}_{m+i}(\mathfrak{X}\setminus S)^{(-1)^{i}q_{n+1-p}(m)\binom{n+1-p}{i}} \cong \\
&\bigotimes_{0\leq p\leq l}\bigotimes_{0\leq i \leq n+1-p}\mathbf{det}(f_{*}\mathcal{L}^{m+i})^{(-1)^{i}q_{n+1-p}(m)\binom{n+1-p}{i}} \cong \\
&\label{sum4} \bigotimes_{0\leq p\leq l}\bigotimes_{0\leq i \leq n+1-p}\bigotimes_{0\leq k \leq n+1}\mathcal{M}_{k}^{(-1)^{i}q_{n+1-p}(m)\binom{n+1-p}{i}\binom{m+i}{k}}\cong\\
&\bigotimes_{0\leq k \leq n+1}\mathcal{M}_{k}^{\sum_{0\leq p\leq l}\sum_{0\leq i \leq n+1-p }(-1)^{i}q_{n+1-p}(m)\binom{n+1-p}{i}\binom{m+i}{k}}.
\end{align}
Next we study the exponent of $\mathcal{M}_{k}$ in (\ref{sum4}).
First we expand the binomial coefficients in of powers of $(m+i)$
 \begin{align*}
 \binom{m+i}{k}=\frac{1}{k!}\sum_{j=0}^{k-1}(-1)^{j}\sigma_{j}(1,2,\dots,k-1)(m+i)^{k-j}.
 \end{align*}
 So that we have
 \begin{align*}
&\sum_{0\leq i \leq n+1-p }(-1)^{i}\binom{n+1-p}{i}\binom{m+i}{k}=\\
&\sum_{0\leq j\leq k-1}\frac{1}{k!}(-1)^{j}\sigma_{j}(1,2,\dots,k-1) \sum_{0\leq i \leq n+1-p}(-1)^{i}\binom{n+1-p}{i}(m+i)^{k-j}=\\
&\sum_{0\leq j\leq k-1}\frac{1}{k!}(-1)^{j+1}\sigma_{j}(1,2,\dots,k-1) P_{n+1-p,k-j}(m)\ .
 \end{align*}
Therefore,
\begin{align*}
&\sum_{0\leq p\leq l}\sum_{0\leq i \leq n+1-p }(-1)^{i}q_{n+1-p}(m)\binom{n+1-p}{i}\binom{m+i}{k}=\\
&\sum_{0\leq p\leq l}\sum_{0\leq j\leq k-1}\frac{1}{k!}(-1)^{j+1}\sigma_{j}(1,2,\dots,k-1) P_{n+1-p,k-j}(m)q_{n+1-p}(m)=\\
&\frac{1}{k!}\sum_{0\leq j\leq k-1}(-1)^{j}\sigma_{j+1}(1,2,\dots,k-1)c_{l.k-j}.
\end{align*}
Which completes the proof of the proposition. 
Recall from the introduction that we have defined $\mathcal{E}_{l}(m)$ as follows 
\begin{align*}
\mathcal{E}_{l}(m):= \sum_{\{0\leq j \leq l\}}\sum_{\{0\leq i \leq n+1-j\}}q_{n+1-j}(m)(-1)^{i}\binom{n+1-j}{i}\mathcal{L}^{m+i} \in \mathbb{Q}[\mathcal{L}].
\end{align*}
(\ref{sum}) says that
\begin{align*}
 &\bigotimes_{0\leq p\leq l}\bigotimes_{0\leq i \leq n+1-p} \mbox{{Hilb}}_{m+i}(\mathfrak{X}\setminus S)^{(-1)^{i}q_{n+1-p}(m)\binom{n+1-p}{i}} \cong \mathbf{det}(\mathcal{E}_{l}(m))\ .
 \end{align*}
This completes the proof of part i) of Theorem 1.

Since the base $S$ is \emph{closed} and $\chi(G)=\{1\}$ ($\chi$ denotes the character group) the action of $G$ on both line bundles must agree, in particular the weights of the respective actions restricted to any one parameter subgroup $\lambda:\mathbb{C}^{*}\rightarrow G$ must agree
\begin{align}\label{wts}
w_{\lambda}(\mathcal{L}_{l})=w_{\lambda}\left(\bigotimes_{0\leq p\leq l}\bigotimes_{0\leq i \leq n+1-p} \mbox{{Hilb}}_{m+i}(\mathfrak{X}\setminus S)^{(-1)^{i}q_{n+1-p}(m)\binom{n+1-p}{i}}\right)\ .
\end{align}
It is easy to see that the weight on the right hand side of (\ref{wts}) is given by $F_{l}(\lambda)$, which is the claim in part ii). Part iii) of theorem 1 is subsumed in the next section.
\begin{center}{\S 3 \sc{Higher Energies}}\end{center}   
Recall from the introduction the \emph{Higher Energy Functionals}.
\begin{align*}
 F_{\om,l}(\vp):=\left(\frac{\sqrt{-1}}{2\pi}\right)^{n}\int_{X}\int_{0}^{1}\frac{\dl}{\dl b}\left(\mbox{Td}(R_{g_{t}}+b\frac{\dl g_{t}}{\dl t}g_{t}^{-1})\mbox{Ch}(F^{\mathcal{E}_{l}(m)}_{H_{t}}+b\frac{\dl H_{t}}{\dl t}H_{t}^{-1})\right)_{|_{b=0}}dt\ .
\end{align*}
The Chern character of $\mathcal{E}_{l}(m)$ is calculated according to the rule
\begin{align}\label{ch}
\begin{split}
&\mbox{Ch}(F^{\mathcal{E}_{l}(m)}_{H_{t}}+b\frac{\dl H_{t}}{\dl t}H_{t}^{-1}):=\\
&\sum_{0\leq j \leq l}\sum_{0\leq i \leq n+1-j}q_{n+1-j}(m)(-1)^{i}\binom{n+1-j}{i}\mbox{Ch}(F^{\mathcal{L}^{m+i}}_{H_{t}}
+b\frac{\dl H_{t}}{\dl t}H_{t}^{-1})\ .
\end{split}
 \end{align}
{The homogeneous decomposition on the right hand side of (\ref{ch}) is in turn given  by the expression}
\begin{align}\label{decomp}
&\mbox{Ch}(F^{\mathcal{L}^{m+i}}_{H_{t}}
+b\frac{\dl H_{t}}{\dl t}H_{t}^{-1})= \sum_{1 \leq k \leq n+1}\mbox{Ch}_{k}(F^{\mathcal{L}^{m+i}}_{H_{t}}
+b\frac{\dl H_{t}}{\dl t}H_{t}^{-1})\ .
\end{align} 
We have dropped the degree zero term  in (\ref{decomp}) as it is killed by $\frac{\dl}{\dl b}$.
In the component of type $(k-1,k-1)$ on the right hand side of (\ref{decomp}) we have defined
$H_{t}= e^{-(m+i)\vp_{t}}h^{(m+i)}$, the induced path of metrics on $\mathcal{L}^{m+i}$.  
$F^{\mathcal{L}^{m+i}}_{H_{t}}$ denotes the curvature.
  \begin{align}
&\frac{\sqrt{-1}}{2\pi} F^{\mathcal{L}^{m+i}}_{H_{t}}=(m+i)(\omega+\frac{\sqrt{-1}}{2\pi}\dl\dlb\varphi_{t})\ .\\
& \dot{H_{t}}H_{t}^{-1}= -(m+i)\dot{\varphi_{t}}\ .
\end{align}
By definition
\begin{align*}
\mbox{Ch}_{k}(F^{\mathcal{L}^{m+i}}_{H_{t}}+b\frac{\dl H_{t}}{\dl t}H_{t}^{-1})=\frac{1}{k!}\mbox{Tr}\left(F^{\mathcal{L}^{m+i}}_{H_{t}}+b\frac{\dl H_{t}}{\dl t}H_{t}^{-1}\right)^{k}.
\end{align*}
Therefore
\begin{align*}
\frac{\dl}{\dl b}\mbox{Ch}_{k}(F^{\mathcal{L}^{m+i}}_{H_{t}}
+b\frac{\dl H_{t}}{\dl t}H_{t}^{-1})_{|_{b=0}}= \frac{(-1)}{(k-1)!}(m+i)^{k}\dot \varphi_{t}\omega^{k-1}_{\varphi_{t}}\ .
\end{align*}
Putting everything together gives
\begin{align*}
&\frac{\dl}{\dl b}\mbox{Ch}(F^{\mathcal{E}_{l}(m)}_{H_{t}}+b\frac{\dl H_{t}}{\dl t}H_{t}^{-1})_{|_{b=0}} \\
&=\sum_{\{0\leq j \leq l\}}\sum_{\{0\leq i \leq n+1-j\}}\sum_{\{1\leq k \leq n+1\}} q_{n+1-j}(m)(-1)^{i}\binom{n+1-j}{i} \frac{(-1)}{(k-1)!}(m+i)^{k}\dot \varphi_{t}\omega^{k-1}_{\varphi_{t}}\\
&=\sum_{\{1\leq k \leq n+1\}}\sum_{\{0\leq j \leq l\}}q_{n+1-j}(m)P_{n+1-j,k}(m)\frac{(-1)}{(k-1)!}\dot \varphi_{t}\omega^{k-1}_{\varphi_{t}}\\
&=\sum_{\{1\leq k \leq n+1\}}\sum_{\{n+1-k \leq j \leq l\}}q_{n+1-j}(m)P_{n+1-j,k}(m)\frac{(-1)}{(k-1)!}\dot \varphi_{t}\omega^{k-1}_{\varphi_{t}}\\
&=\sum_{\{n+1-l \leq k \leq n+1\}}\frac{(-1)}{(k-1)!}c_{l,k}\dot \varphi_{t}\omega^{k-1}_{\varphi_{t}}\ .
\end{align*}
Therefore
\begin{align}
\frac{\dl}{\dl b}\mbox{Ch}_{k}(F^{\mathcal{E}_{l}(m)}_{H_{t}}+b\frac{\dl H_{t}}{\dl t}H_{t}^{-1})_{|_{b=0}} = \frac{(-1)}{(k-1)!}c_{l,k}\dot \varphi_{t}\omega^{k-1}_{\varphi_{t}} \ .
\end{align}

In the same way one shows that
\begin{align}
\mbox{Ch}(F^{\mathcal{E}_{l}(m)}_{H_{t}}) =
 \sum _{\{n+1-l \leq k \leq n\}}\frac{c_{l,k}}{k!}\omega_{ \varphi_{t}}^{k}.
\end{align}
The top dimensional component of  the transgressed GRR integrand is given by     
 \begin{align*}
&\frac{\dl}{\dl b}\mbox{Ch}(F^{\mathcal{E}_{l}(m)}_{H_{t}}+b\frac{\dl H_{t}}{\dl t}H_{t}^{-1})_{|_{b=0}}\mbox{Td}(R_{g_{t}})+
\frac{\dl}{\dl b}\mbox{Td}(R_{g_{t}}+b\frac{\dl g_{t}}{\dl t}g_{t}^{-1})_{|_{b=0}}\mbox{Ch}(F^{\mathcal{E}_{l}(m)}_{H_{t}})\ =\\
&\sum_{\{n+1-l\leq k\leq n+1\}}\mbox{Td}_{n+1-k}(R_{g_{t}})\frac{\dl}{\dl b}\mbox{Ch}_{k}(F^{\mathcal{E}_{l}(m)}_{H_{t}}+b\frac{\dl H_{t}}{\dl t}H_{t}^{-1})_{|_{b=0}}\\
&+\sum_{\{n+1-l\leq k \leq n\}}\frac{\dl}{\dl b}\mbox{Td}_{n+1-k}(R_{g_{t}}+b\frac{\dl g_{t}}{\dl t}g_{t}^{-1})_{|_{b=0}}\mbox{Ch}_{k}(F^{\mathcal{E}_{l}(m)}_{H_{t}})\ 
\end{align*}
We summarize what we have in the following lemma.
 \begin{proposition}
\begin{align*} 
\begin{split}
&F_{\omega,l}(\varphi)=
\int_{X}{\int_{0}}^{1}\sum_{n+1-l\leq k\leq n}\frac{c_{l,k}}{k!}\frac{\dl}{\dl b}Td_{n+1-k}(\frac{\sqrt{-1}}{2\pi}R_{g_{t}}+b\frac{\dl g_{t}}{\dl t}g_{t}^{-1})|_{b=0}
{\omega_{t}}^{k}dt\\
&-\int_{X}{\int_{0}}^{1} \sum_{n+1-l\leq k \leq n+1}\frac{k c_{l,k} }{k!}Td_{n+1-k}(\frac{\sqrt{-1}}{2\pi}R_{g_{t}})\frac{\dl}{\dl t}\varphi_{t} {\omega_{t}}^{k-1} dt\ .
\end{split}
\end{align*}
 \end{proposition}
In particular when  $l=2$  we have the following.
\begin{align}\label{efl}
\begin{split}
&\frac{ \dl F_{\om,2}(\vp_{t})}{\dl t}=\int_{X}\Bigg( \frac{c_{l,n-1}}{(n-1)!}\frac{\dl}{\dl b}Td_{2}(\frac{\sqrt{-1}}{2\pi}R_{\varphi}+b(\dot\vp_{i\overline{j}})g_{\varphi}^{-1})|_{b=0}{\omega^{n-1}_{\varphi}}-\\
&  \frac{c_{l,n-1} }{(n-2)!}Td_{2}(\frac{\sqrt{-1}}{2\pi}R_{\varphi})\dot\vp {\omega^{n-2}_{\varphi}}
 -\frac{c_{l,n} }{(n-1)!}Td_{1}(\frac{\sqrt{-1}}{2\pi}R_{\varphi})\dot\vp {\omega^{n-1}_{\varphi}} - \frac{ c_{l,n+1} }{n!} \dot\vp {\omega^{n}_{\varphi}}\Bigg)
\end{split}
\end{align}
In order to compare $F_{\om,2}$ with $E_{\om,1}$ we need to make the right hand side of (\ref{efl}) more explicit.
\begin{proposition}
\begin{align}
 & i) \ \frac{\dl}{\dl b}Td_{2}(\frac{\sqrt{-1}}{2\pi}R_{\varphi}+b(\dot\vp_{i\overline{j}})g_{\varphi}^{-1})|_{b=0}{\omega_{\varphi}}^{n-1}=\frac{1}{6n}\Delta_{\varphi_{t}}\dot\varphi_{t}\mbox{\emph{Scal}}(\varphi_{t})\omega_{\vp_{t}}^{n}\\
&ii) \ Td_{2}(\frac{\sqrt{-1}}{2\pi}R_{\varphi})\dot \vp {\omega_{\varphi}}^{n-2}=\frac{1}{12n(n-1)}\dot\varphi_{t}\left(3\mbox{\emph{Scal}}(\vp_{t})^{2}-4||\mbox{\emph{Ric}}_{\vp}||^{2}+||R_{\vp}||^{2}\right)\ .
\end{align}
\end{proposition}
For a vector bundle $\mathcal{E}\ra X$ recall that the second Todd class is given by the formula
\begin{align*}
12 \mbox{Td}_{2}(\mathcal{E})=  \frac{3}{2}c_{1}(\mathcal{E})^{2}-\mbox{Ch}_{2}(\mathcal{E}).\\
\end{align*}
On a K\"ahler manifold the Chern-Weil theory  implies that for $\mathcal{E}=T^{1,0}_{X}$ we get
\begin{align}\label{td2}
12Td_{2}(\frac{\sqrt{-1}}{2\pi}R_{\varphi})=\frac{3}{2}\mbox{\mbox{Ric}}(\om_{\vp})^{2}- \frac{1}{2}\mbox{{Tr}}(\frac{\sqrt{-1}}{2\pi}R_{\varphi})^{2}.
\end{align}
 Therefore we have
 \begin{align*}
&12\int_{X}\frac{\dl}{\dl b}\mbox{Td}_{2}\left(\frac{\sqrt{-1}}{2\pi}R_{\varphi}+b(\dot\vp_{i\overline{j}})g_{\varphi}^{-1}\right)|_{b=0}\ \om_{\vp}^{n-1} =\\
&\frac{3}{2}\int_{X}\frac{\dl}{\dl b}\Bigg\{\mbox{Tr}\left(\frac{\sqrt{-1}}{2\pi}R_{\varphi}+b(\dot\vp_{i\overline{j}})g_{\varphi}^{-1}\right)\Bigg\}^{2}|_{b=0}\ \om_{\vp}^{n-1}+\\
& -\frac{1}{2}\int_{X}\frac{\dl}{\dl b}\mbox{Tr}\Bigg\{\frac{\sqrt{-1}}{2\pi}R_{\varphi}+b(\dot\vp_{i\overline{j}})g_{\varphi}^{-1}\Bigg\}^{2}|_{b=0}\ \om_{\vp}^{n-1} =\\
&3\int_{X}\mbox{Tr}\left(\frac{\sqrt{-1}}{2\pi}R_{\varphi}\right)\Delta_{\vp}\dot\vp\ \om_{\vp}^{n-1}-\int_{X}\mbox{Tr}\left(\frac{\sqrt{-1}}{2\pi}R_{\varphi}(\dot\vp_{i\overline{j}})g^{-1}_{\vp}\right)\om_{\vp}^{n-1}=\\
& 3\int_{X}\mbox{Ric}(\omega_{\varphi})\Delta_{\varphi}\dot\vp\ \om_{\vp}^{n-1}-
\int_{X}\mbox{Tr}(\frac{\sqrt{-1}}{2\pi}R_{\varphi}(\dot\vp_{i\overline{j}})g_{\varphi}^{-1})\ \om_{\vp}^{n-1}=\\
&\frac{3}{n}\int_{X}\dot \vp\Delta_{\vp}\mbox{Scal}_{\vp}{\om^{n}}_{\vp}-\int_{X}\mbox{Tr}(\frac{\sqrt{-1}}{2\pi}R_{\varphi}(\dot\vp_{i\overline{j}})g_{\varphi}^{-1})\ \om_{\vp}^{n-1}\ .
\end{align*}
Recall the following fact.  
 \begin{align}\label{contract}
\Theta\wedge\omega^{n-k}=\frac{(n-k)!}{n!}\Lambda^{k}\Theta \ \omega^{n} \ .
\end{align}
$\Lambda$ denotes contraction against the K\"ahler form $\om$ and $\Theta$ denotes a $(k,k)$ form on $X$. 
An application of (\ref{contract}) gives that 
 \begin{align}\label{trace}
\int_{X}n\mbox{Tr}(\frac{\sqrt{-1}}{2\pi}R_{\varphi}(\dot\vp_{i\overline{j}})g_{\varphi}^{-1})\omega_{\varphi}^{n-1}=  \int_{X}\dot\vp\sum_{1\leq i,j,k \leq n}\frac{\dl^{2}}{\dl z_{j}\dl \overline{z}_{k}}\left(\Lambda R_{g}(i,j)g^{k\overline{i}}\mbox{det}(g)\right)\frac{\omega^{n}_{g}}{\mbox{det}(g)}
 \end{align}
On the right hand side of (\ref{trace}) we have temporarily replaced $g_{\vp}$ by $g$.

Next we need some local formulas from K\"ahler geometry which we collect below.
\begin{align*}
&R_{g}(i,j):=\dlb \{\dl(g) g^{-1}\}(i,j) =\sum_{k,l}R^{k\overline{l}}_{i\overline{j}}\ dz_{k}\wedge d\overline{z}_{l}\quad (\mbox{full curvature tensor})\\
&\mbox{Ric}_{g}:=\mbox{Tr}(R_{g})=\sum_{\{1\leq k,l \leq n\}}\sum_{\{1\leq i\leq n\}}R^{k\overline{l}}_{i\overline{i}}dz_{k}\wedge d\overline{z}_{l}\quad (\mbox{Ricci curvature})\\
&\mbox{Scal}_{\omega}:= \mbox{Tr}(\Lambda R_{g})=\sum_{1\leq i,k,l\leq n}g^{l\overline{k}}R^{k\overline{l}}_{i\overline{i}}\ \quad (\mbox{scalar curvature})\ .\\
\end{align*}
 
At the center $o$ of a normal coordinate system we have that
\begin{align*}
&R^{k\overline{l}}_{i\overline{j}}(o)=-\frac{\dl^{2}  g_{i\overline{j}}}{\dl z_{k}\dl \overline{z}_{l}}(o)=\frac{\dl^{2}  g^{i\overline{j}}}{\dl z_{k}\dl \overline{z}_{l}}(o)\\
 &\Delta_{g}\mbox{Scal}_{\omega}(o)=\sum_{1\leq i,j,k\leq n}\frac{\dl^{2}R^{k\overline{k}}_{i\overline{i}}}{\dl z_{j}\dl \overline{z}_{j}}(o)
+\sum_{1\leq i,j,k,l\leq n}R^{k\overline{l}}_{i\overline{i}}(o)R^{l\overline{k}}_{j\overline{j}}(o)\ .
\end{align*}
Since $\mbox{Ric}_{g}$ is closed we have the identity
\begin{align}
\sum_{1\leq i \leq n}\frac{\dl^{2}R^{k\overline{l}}_{i\overline{i}}}{\dl z_{p}\dl \overline{z}_{q}}=\sum_{1\leq i \leq n}\frac{\dl^{2}R^{p\overline{q}}_{i\overline{i}}}{\dl z_{k}\dl \overline{z}_{l}}\ .
\end{align}\label{switch}
The commutation $ R_{i\overline{j}k\overline{l}}=R_{k\overline{l}i\overline{j}}$ implies that at $o$ we have
\begin{align}\label{updown}
\frac{\dl^{2}R^{k\overline{l}}_{i\overline{j}}}{\dl z_{p}\dl \overline{z}_{q}}=\frac{\dl^{2}R^{i\overline{j}}_{k\overline{l}}}{\dl z_{p}\dl \overline{z}_{q}}+\sum_{1\leq m \leq n}(R^{k\overline{l}}_{i\overline{m}}R^{p\overline{q}}_{m\overline{j}}-R^{i\overline{j}}_{k\overline{m}}R^{p\overline{q}}_{m\overline{l}})\ .
\end{align} 
 Therefore we get
 \begin{align*}
&\frac{\dl^{2}}{\dl z_{j}\dl \overline{z}_{k}}\left(\Lambda R_{g}(i,j)g^{k\overline{i}}\mbox{det}(g)\right)=\\
&\frac{\dl^{2}R^{p\overline{p}}_{i\overline{j}}}{\dl z_{j}\dl \overline{z}_{i}}+R^{q\overline{p}}_{i\overline{j}}R^{j\overline{i}}_{p\overline{q}}
+R^{k\overline{i}}_{j\overline{k}}R^{p\overline{p}}_{i\overline{j}}- R^{p\overline{p}}_{i\overline{j}}R^{j\overline{i}}_{l\overline{l}}=\\
&\frac{\dl^{2}R^{i\overline{j}}_{p\overline{p}}}{\dl z_{j}\dl \overline{z}_{i}}+R^{p\overline{p}}_{i\overline{m}}R^{j\overline{i}}_{m\overline{j}}- R^{i\overline{j}}_{p\overline{m}}R^{j\overline{i}}_{m\overline{p}}+R^{q\overline{p}}_{i\overline{j}}R^{j\overline{i}}_{p\overline{q}}
+R^{k\overline{i}}_{j\overline{k}}R^{p\overline{p}}_{i\overline{j}}- R^{p\overline{p}}_{i\overline{j}}R^{j\overline{i}}_{l\overline{l}}\ .
\end{align*}
Next use the symmetries
\begin{align*}
&\frac{\dl^{2}R^{i\overline{j}}_{p\overline{p}}}{\dl z_{j}\dl \overline{z}_{i}}=\frac{\dl^{2}R^{j\overline{j}}_{p\overline{p}}}{\dl z_{i}\dl \overline{z}_{i}}\\
&R^{i\overline{j}}_{k\overline{l}}=R^{k\overline{j}}_{i\overline{l}}=R^{i\overline{l}}_{k\overline{j}}=R^{k\overline{l}}_{i\overline{j}}\ .
\end{align*}
 To deduce that
\begin{align*}
&\frac{\dl^{2}}{\dl z_{j}\dl \overline{z}_{k}}\left(\Lambda R_{g}(i,j)g^{k\overline{i}}\mbox{det}(g)\right)=\\
&\frac{\dl^{2}R^{j\overline{j}}_{p\overline{p}}}{\dl z_{i}\dl \overline{z}_{i}}+R^{p\overline{p}}_{i\overline{m}}R^{j\overline{j}}_{m\overline{i}}
- R^{i\overline{j}}_{p\overline{m}}R^{j\overline{i}}_{m\overline{p}}+R^{q\overline{p}}_{i\overline{j}}R^{p\overline{q}}_{j\overline{i}}
+R^{k\overline{k}}_{j\overline{i}}R^{p\overline{p}}_{i\overline{j}}- R^{p\overline{p}}_{i\overline{j}}R^{l\overline{l}}_{j\overline{i}}=\\
&\frac{\dl^{2}R^{j\overline{j}}_{p\overline{p}}}{\dl z_{i}\dl \overline{z}_{i}}+R^{k\overline{k}}_{j\overline{i}}R^{p\overline{p}}_{i\overline{j}}=\Delta_{g}\mbox{Scal}_{\omega}(o)\ .
\end{align*}
So that we have
\begin{align*}
\int_{X}\mbox{Tr}(\frac{\sqrt{-1}}{2\pi}R_{\varphi}(\dot\vp_{i\overline{j}})g_{\varphi}^{-1})\omega_{\varphi}^{n}=\frac{1}{n}\int_{X}\dot\vp\Delta_{\varphi}\mbox{Scal}_{\varphi}\omega_{\varphi}^{n}\ .
\end{align*}
This establishes $i)$.

 By definition we have
 \begin{align}\label{ricci}
 &{\mbox{Ric}_{g}}^{2}=R^{k\overline{l}}_{i\overline{i}}R^{p\overline{q}}_{j\overline{j}} \ dz_{k}\wedge d\overline{z}_{l}\wedge dz_{p}\wedge d\overline{z}_{q}\\
 & \label{trace2} \mbox{Tr}(\mbox{R}_{g}^{2})=R^{k\overline{l}}_{i\overline{j}}R^{p\overline{q}}_{j\overline{i}}\ dz_{k}\wedge d\overline{z}_{l}\wedge dz_{p}\wedge d\overline{z}_{q} \ .
 \end{align}
 It is easy to see that
\begin{align}\label{contract2}
\Lambda^{2}\{dz_{k}\wedge d\overline{z}_{l}\wedge dz_{p}\wedge d\overline{z}_{q}\}=2\{\delta_{kl}\delta_{pq}-\delta_{kq}\delta_{pl}\}\ .
\end{align}
 Applying (\ref{contract2}) to (\ref{ricci}) and (\ref{trace2}) respectively gives
 \begin{align}\label{joanna}
\begin{split}
&{\Lambda_{g}}^{2}{\mbox{Ric}_{g}}^{2}=2(S_{g}^{2}- ||\mbox{Ric}_{g}||^{2})\\
\ \\
&{\Lambda_{g}}^{2}\mbox{Tr}({R_{g}}^{2})= 2(||\mbox{Ric}_{g}||^{2}-||R_{g}||^{2})\ . 
\end{split}
\end{align}
Recall the identity
\begin{align} \label{joanna2}
12\mbox{Td}_{2}(\frac{\sqrt{-1}}{2\pi}R_{\varphi})\om_{\vp}^{n-2}&=\frac{3}{2}\mbox{\mbox{Ric}}(\om_{\vp})^{2}\om_{\vp}^{n-2}- \frac{1}{2}\mbox{{Tr}}(\frac{\sqrt{-1}}{2\pi}R_{\varphi})^{2}\om_{\vp}^{n-2} \ .
\end{align}
An application of (\ref{contract}) and (\ref{joanna}) to the right hand side of (\ref{joanna2}) shows that
\begin{align*}
&12\mbox{Td}_{2}(\frac{\sqrt{-1}}{2\pi}R_{\varphi})\om_{\vp}^{n-2}=\\
&\frac{3}{n(n-1)}(S_{g}^{2}- ||\mbox{Ric}_{g}||^{2})
- \frac{1}{n(n-1)}(||\mbox{Ric}_{g}||^{2}-||R_{g}||^{2})\\ 
&=\frac{1}{n(n-1)}\left(3\mbox{Scal}_{\vp}^{2}-4||\mbox{Ric}_{\vp}||^{2}+||R_{\vp}||^{2}\right)\ .
\end{align*}
We leave the computation of the coefficients $c_{2,j}\ (n-1\leq j \leq n+1)$ to the reader. This establishes proposition 2.

\begin{center}{\S 4 \sc{Higher Energy Asymptotics}}\end{center}
 In this section we exhibit the energy $F_{2,\om}$ as a singular Hermitian metric on the sheaf $\mathcal{L}_{2}$. This is a carried out using the method of Tian \emph{mutatis mutandis}.
Let $f^{-1}(o)=X_{o}\subset \cpn$, where $o\in {S}_{\infty}:=S\setminus \Delta$, where $\Delta$ denotes the discriminant locus of the family.
We define for any $z\in S$
\begin{align*}
GX_{z}:= \{(\sigma,y)\in G \times \cpn:y\in \sigma X_{z}\} \ . 
\end{align*}
Then we have the following diagram, where $p_{z}$ denotes the evaluation map, i.e.
$p_{z}(\sigma):= \sigma z$.

 \begin{diagram}
p_{z}^{*}(\mathfrak{X})\cong GX_{z}&\rTo^{p_{z,2}}&\mathfrak{X}&\rInto^{\iota}&\mathbb{P}(f_{*}{L}^{\vee})\cong B\times \cpn&\rTo^{p_{2}}&\cpn \\
\dTo ^{p_{z,1}}&&\dTo^{f}&\ldTo_{\pi}&\\
 G&\rTo ^{p_{z}}&S&
\end{diagram}
Given $z\in S\setminus \Delta$ we can consider $T^{1,0}_{X_{z}}$, the holomorphic tangent bundle of the fiber $X_{z}$. These fit together holomorphically into a vector bundle $\mathcal{V}$ on $\mathfrak{X}\setminus f^{-1}(\Delta)$.
We have the following exact complex (denoted by $(A^{\bull}\ ,\ \dl_{\bull})$ ) over $\mathfrak{X}\setminus f^{-1}(\Delta)$.
 
\[ \begin{CD} 
0@>>>\mathcal{V}@>\iota_{*}>>T^{1,0}_{\mathfrak{X}}|_{\mathfrak{X}\setminus f^{-1}(\Delta)}@>f_{*}>>f^{*}T^{1,0}_{S}|_{\mathfrak{X}\setminus f^{-1}(\Delta)}@>>>0 \ . \end{CD}
\]
Let $\om_{\mathfrak{X}}$ and $\om_{S}$ denote any K\"ahler metrics on $\mathfrak{X}$ and $S$ respectively. The Fubini study form $\om$ on $\cpn$ induces a K\"ahler form on each of the smooth fibers of $f$ and so induces a Hermitian metric  (which we denote by $h_{\mathcal{V}}$ ) on $\mathcal{V}$.

More generally let
\[ \begin{CD} 
0@>>>E^{0}@>v_{0}>>E^{1}@>v_{1}>> \dots @>>> E^{i}@>v_{i}>> E^{i+1}@>>> \dots @>v_{l-1}>> E^{l}@>>>0
 \end{CD}
\]
be a complex of holomorphic \emph{Hermitian} vector bundles on a complex manifold $B$. The metric on $E^{i}$ will be denoted by $h_{i}$, the corresponding holomorphic Hermitian connection by $\nabla_{i}$.  $\nabla^{2}_{i}$ is the curvature.
Let $N$ denote the \emph{number operator} of the complex, i.e. $N$ acts by multiplication by $j$ on $E^{j}\ (0\leq j \leq l)$. 
For $u\geq 0$, let $A_{u}$ be the Quillen superconnection
\begin{align*}
A_{u}:= \nabla+\sqrt{u}V\ , \quad V:= v+v^{*}\ .
\end{align*}
\begin{definition} (Bismut, Gillet, Soul\'e \cite{bgs1})\\
For $s\in \mathbb{C}$, $Re(s)>0$, let $\zeta_{E^{\bull}}(s)$ be the collection of forms on B defined by
\begin{align*}
\zeta_{E^{\bull}}(s):= \frac{-1}{\Gamma(s)}\int_{0}^{+\infty}u^{s-1}\mbox{Tr}_{s}[N\mbox{exp}(-A^{2}_{u})]du
\end{align*}
$\zeta_{E^{\bull}}(s)$ is actually holomorphic on all of $\mathbb{C}$.
\end{definition}
 
\noindent \textbf{Theorem} \emph{(Bismut, Gillet, Soul\'e \cite{bgs1} Theorem 1.15)}\\ 
\emph{ On the base $B$ there is a pointwise identity of forms provided that the complex $(E^{\bull}\ ,\ v_{\bull})$ is acyclic.}
\begin{align}\label{bgs}
\sum_{0\leq j \leq l}(-1)^{j}\mbox{Ch}(E^{j},\ h_{j})=  
\dl\dlb{\dot \zeta}_{E^{\bull}}(0)\ .
\end{align}
${\dot \zeta}_{E^{\bull}}(0)^{\{p,p\}}$ denotes  the component of ${\dot \zeta}_{E^{\bull}}(0)$ of degree $(p,p)$.
 Direct application of (\ref{bgs}) to the adjunction complex gives forms ${\dot \zeta}_{A^{\bull}}(0)^{\{p,p\}}$  \emph{defined along } $\mathfrak{X}\setminus f^{-1}(\Delta)$ satisfying
 \begin{align}
\begin{split}\label{appl}
&\dl\dlb {\dot \zeta}_{A^{\bull}}(0)^{\{p,p\}}=\\
 &\mbox{Ch}(\mathcal{V},\om)^{\{p+1,p+1\}}-\mbox{Ch}(T^{1,0}_{\mathfrak{X}},\om_{\mathfrak{X}})^{\{p+1,p+1\}}+f^{*}\mbox{Ch}(T^{1,0}_{S},\om_{S})^{\{p+1,p+1\}}\ .
 \end{split}
 \end{align}
When $p=0$ we can describe the function ${\dot \zeta}_{A^{\bull}}(0)^{\{0,0\}}$ in the following way.

Given $z\in S\setminus \Delta$ we can consider $K_{X_{z}}$, the canonical bundle of the fiber $X_{z}$. These fit together holomorphically into a line bundle $K_{\infty}\ (\cong \mathbf{det}(\mathcal{V})^{-1})$ on $\mathfrak{X}\setminus f^{-1}(\Delta)$.
On the other hand the relative canonical bundle  $K_{f}$ of the map $f$ is given by 
 \begin{align*}
{K}_{f}:= K_{\mathfrak{X}}\otimes f^{*} K_{S}^{-1} \ .
\end{align*}
When we restrict  ${K}_{f}$ to $\mathfrak{X}\setminus f^{-1}(\Delta)$ we have an isomorphism
\begin{align*}
{K}_{f}  \cong K_{\infty}\ .
\end{align*}
$\iota^{*}p_{2}^{*}\omega_{FS}$ restricts to a K\"ahler  metric on $f^{-1}(z)$ $(z \in {S}_{\infty})$ and hence induces a Hermitian metric on the bundle $K_{\infty}$.  We denote 
 its curvature
 by $R( \iota^{*}p_{2}^{*} (\omega_{FS}))$.
 The K\"ahler metrics on
$\mathfrak{X}$ and $S$ 
 induce a metric on the relative canonical bundle $K_{f}$.
We let $R_{f}$ denote its curvature
\begin{align*}
R_{\mathfrak{X}/ S}:= R(\om_{\mathfrak{X}})-f^{*}R(\om_{S}).
\end{align*}
 In this way we obtain \emph{two} metrics on the relative canonical bundle over the smooth locus. The curvatures of these metrics are \emph{not} the same.
The relation between them is given in the following ``\emph {$\dl\dlb$ lemma along the fibers}''.\\
 
\emph{There is a smooth function $\Psi:\mathfrak{X}\setminus f^{-1}(\Delta)\rightarrow \mathbb{R}$ such that}
\begin{align}
\begin{split}
& R_{\mathfrak{X}/ S}  + \frac{\sqrt{-1}}{2\pi}\dl\dlb \Psi = R( \iota^{*}p_{2}^{*} (\omega_{FS}))\ . 
 \end{split}
\end{align}
Then (up to additive constants) we have the identity
\begin{align}
-{\dot \zeta}_{A^{\bull}}(0)^{\{0,0\}}=\Psi \ .
\end{align}
\begin{lemma}
There is a continuous metric $|| \  ||_{2}$ on the line bundle $\mathcal{L}_{2}$ such that in the sense of currents we have
\begin{align}
\begin{split}
&\frac{\sqrt{-1}}{2\pi}\dl\dlb \mbox{\emph{log}}(|| \  ||^{2}_{2})= \\
&f_{*}\Bigg\{\frac{c_{2,n-1}}{(n-1)!}\Bigg(\frac{1}{8}(R_{\mathfrak{X}/ S})^{2}p_{2}^{*}\om_{FS}^{n-1} -\frac{1}{12}(\mbox{\emph{Ch}}_{2}(T^{1,0}_{\mathfrak{X}},\om_{\mathfrak{X}})-f^{*}\mbox{\emph{Ch}}_{2}(T^{1,0}_{S},\om_{S}))p_{2}^{*}\om_{FS}^{n-1}\Bigg)\\
&- \frac{c_{2,n}}{n!} (R_{\mathfrak{X}/ S})p_{2}^{*}\om_{FS}^{n} +\frac{c_{2,n+1}}{(n+1)!}p_{2}^{*}\om_{FS}^{n+1}
\Bigg\} \ .
\end{split}
\end{align}
\end{lemma}
Next we define a Hermitian metric $h_{G}$ on the vector bundle $\pi_{2}^{*}T^{1,0}_{X_{z}}$ over the product $G\times X_{z}$
\begin{align}
h_{G}|_{\sigma \times X_{z}}:=( \omega_{FS}+\frac{\sqrt{-1}}{2\pi}\dl\dlb \varphi_{\sigma})|_{X_{z}}\ .
\end{align}
Then we have the corresponding Chern-Weil forms
\begin{align}
\mbox{Ch}_{k}(\pi_{2}^{*}T^{1,0}_{X_{z}}, h_{G})\ .
\end{align}
\begin{lemma}\label{ddbar}
Let $\psi$ be a compactly supported form on $G$. Then
\begin{align}
\begin{split}
&\int_{G}F_{2,\om}(\vp_{\sigma})\dl\dlb \psi = \\
&\int_{G\times X_{z}}\Bigg\{\frac{c_{2,n-1}}{(n-1)!}\Bigg(\frac{1}{8}\mbox{c}_{1}(\pi_{2}^{*}T^{1,0}_{X_{z}}, h_{G})   ^{2}\pi_{2}^{*}\om_{FS}^{n-1} -
\frac{1}{12}\mbox{\emph{Ch}}_{2}(\pi_{2}^{*}T^{1,0}_{X_{z}}, h_{G})\pi_{2}^{*}\om_{FS}^{n-1}\Bigg)\\
&- \frac{c_{2,n}}{n!}\mbox{c}_{1}(\pi_{2}^{*}T^{1,0}_{X_{z}}, h_{G})\pi_{2}^{*}\om_{FS}^{n}
 +\frac{c_{2,n+1}}{(n+1)!}\pi_{2}^{*}\om_{FS}^{n+1}
\Bigg\}\wedge \pi_{1}^{*}(\psi)
\end{split}
\end{align}
\end{lemma}
Next observe that $GX_{z}$ is biholomorphic to $G\times X_{z}$. It follows that we have an identity
\begin{align*}
p^{*}_{z,2}\mbox{Ch}_{k}(\mathcal{V}, h_{\mathcal{V}})= \mbox{Ch}_{k}(\pi_{2}^{*}T^{1,0}_{X_{z}}, h_{G})\ .
\end{align*}
 \noindent It then follows from (\ref{appl}) that we have
\begin{align}\label{fiber}
\begin{split}
& \mbox{c}_{1}(\pi_{2}^{*}T^{1,0}_{X_{z}}, h_{G})  =p^{*}_{z,2}\Bigg(-R_{\mathfrak{X}/ S}+\frac{\sqrt{-1}}{2\pi}\dl\dlb {\dot \zeta}_{A^{\bull}}(0)^{\{0,0\}}\Bigg)\\
&  \mbox{Ch}_{2}(\pi_{2}^{*}T^{1,0}_{X_{z}}, h_{G})= p^{*}_{z,2}\Bigg(\mbox{Ch}_{2}(T^{1,0}_{\mathfrak{X}},\om_{\mathfrak{X}})-f^{*}\mbox{{Ch}}_{2}(T^{1,0}_{S},\om_{S})+\dl\dlb {\dot \zeta}_{A^{\bull}}(0)^{\{1,1\}}\Bigg) \ .
\end{split}
\end{align}
This implies the 
\begin{corollary}
The function
\begin{align}\label{pluri}
\sigma \in G \ra D(\sigma):= F_{2,\om|_{X_{z}}}(\sigma)-\mbox{\emph{log}}\Bigg(e^{\Psi_{S}(\sigma z)}\frac{||\ ||_{2}^{2}(\sigma z)}{|| \ ||^{2}_{2}(z)}\Bigg)
\end{align}
is pluriharmonic.   $\Psi_{S}:S\setminus \Delta \ra \mathbb{R}$ is given by
\begin{align}
\begin{split}
&\Psi_{S}(z):=\kappa_{1}\int_{X_{z}}\Bigg(-\frac{\sqrt{-1}}{\pi}{\dot \zeta}_{A^{\bull}}(0)^{\{0,0\}}R_{\mathfrak{X}/ S} -\frac{1}{4\pi^{2}} {\dot \zeta}_{A^{\bull}}(0)^{\{0,0\}}\dl\dlb{\dot \zeta}_{A^{\bull}}(0)^{\{0,0\}}\Bigg)\pi_{2}^{*}\om^{n-1}\\
&+\kappa_{2}\int_{X_{z}}{\dot \zeta}_{A^{\bull}}(0)^{\{1,1\}}\pi_{2}^{*}\om^{n-1}+\kappa_{3}\int_{X_{z}}\frac{\sqrt{-1}}{2\pi}\dl\dlb {\dot \zeta}_{A^{\bull}}(0)^{\{0,0\}}\pi_{2}^{*}\om^{n}.\\
& \ \\
&8(n-1)!\kappa_{1}=c_{2,n-1} \quad -12(n-1)!\kappa_{2}=c_{2,n-1}\quad -n!\kappa_{3}=c_{2,n}.
\end{split}
\end{align}
\end{corollary}
Standard arguments  show that the right hand side of (\ref{pluri}) is identically \emph{zero}. This exhibits the new energies as singular metrics on the lines $\mathcal{L}_{l}\ 0\leq l \leq n+1$, although we have only carried this out in detail when $l=2$. In order to conclude the proof of Theorem 2 let $\sigma=\lambda(t)$ on the right hand side of (\ref{pluri}) and let $t\ra 0$. Since we have assumed that the limit cycle $X_{z}^{\lambda(0)}$ is $C^{\infty}$ the term $\Psi$ is bounded. \newline
\emph{Q.E.D.}

\bibliography{ref}

\begin{thebibliography}{10}

\bibitem{cay}
A.Cayley.
\newblock On the theory of elimination.
\newblock {\em Cambridge and Dublin Math Journal}, 3, 1848.

\bibitem{bandmab}
Shigetoshi Bando and Toshiki Mabuchi.
\newblock Uniqueness of {E}instein {K}\"ahler metrics modulo connected group
  actions.
\newblock In {\em Algebraic geometry, Sendai, 1985}, volume~10 of {\em Adv.
  Stud. Pure Math.}, pages 11--40. North-Holland, Amsterdam, 1987.

\bibitem{bgs1}
J.-M. Bismut, H.~Gillet, and C.~Soul{\'e}.
\newblock Analytic torsion and holomorphic determinant bundles. {I}.
  {B}ott-{C}hern forms and analytic torsion.
\newblock {\em Comm. Math. Phys.}, 115(1):49--78, 1988.

\bibitem{cat}
David Catlin.
\newblock The {B}ergman kernel and a theorem of {T}ian.
\newblock In {\em Analysis and geometry in several complex variables (Katata,
  1997)}, Trends Math., pages 1--23. Birkh\"auser Boston, Boston, MA, 1999.

\bibitem{chentianricci}
X.~X. Chen and G.~Tian.
\newblock Ricci flow on {K}\"ahler-{E}instein manifolds.
\newblock {\em Duke Math. J.}, 131, 2006.

\bibitem{dingtian}
Wei~Yue Ding and Gang Tian.
\newblock K\"ahler-{E}instein metrics and the generalized {F}utaki invariant.
\newblock {\em Invent. Math.}, 110(2):315--335, 1992.

\bibitem{skdproj}
S.~K. Donaldson.
\newblock Scalar curvature and projective embeddings. {I}.
\newblock {\em J. Differential Geom.}, 59(3):479--522, 2001.

\bibitem{skdtoric}
S.~K. Donaldson.
\newblock Scalar curvature and stability of toric varieties.
\newblock {\em J. Differential Geom.}, 62(2):289--349, 2002.

\bibitem{fogarty}
John Fogarty.
\newblock Truncated {H}ilbert functors.
\newblock {\em J. Reine Angew. Math.}, 234:65--88, 1969.

\bibitem{tableaux}
William Fulton.
\newblock {\em Young tableaux}, volume~35 of {\em London Mathematical Society
  Student Texts}.
\newblock Cambridge University Press, Cambridge, 1997.
\newblock With applications to representation theory and geometry.

\bibitem{fulton}
William Fulton.
\newblock {\em Intersection theory}, volume~2 of {\em Ergebnisse der Mathematik
  und ihrer Grenzgebiete. 3. Folge. A Series of Modern Surveys in Mathematics
  [Results in Mathematics and Related Areas. 3rd Series. A Series of Modern
  Surveys in Mathematics]}.
\newblock Springer-Verlag, Berlin, second edition, 1998.

\bibitem{globmod}
D.~Gieseker.
\newblock Global moduli for surfaces of general type.
\newblock {\em Invent. Math.}, 43(3):233--282, 1977.

\bibitem{ksz}
M.~M. Kapranov, B.~Sturmfels, and A.~V. Zelevinsky.
\newblock Chow polytopes and general resultants.
\newblock {\em Duke Math. J.}, 67(1):189--218, 1992.

\bibitem{detdiv}
Finn~Faye Knudsen and David Mumford.
\newblock The projectivity of the moduli space of stable curves {I}:
  {P}reliminaries on ``det'' and ``{D}iv''.
\newblock {\em Math. Scand.}, 39(1), 1976.

\bibitem{mabuchi}
Toshiki Mabuchi.
\newblock {$K$}-energy maps integrating {F}utaki invariants.
\newblock {\em Tohoku Math. J. (2)}, 38(4):575--593, 1986.

\bibitem{git}
D.~Mumford, J.~Fogarty, and F.~Kirwan.
\newblock {\em Geometric invariant theory}, volume~34 of {\em Ergebnisse der
  Mathematik und ihrer Grenzgebiete (2) [Results in Mathematics and Related
  Areas (2)]}.
\newblock Springer-Verlag, Berlin, third edition, 1994.

\bibitem{sopv}
David Mumford.
\newblock Stability of projective varieties.
\newblock {\em Enseignement Math. (2)}, 23(1-2):39--110, 1977.

\bibitem{cms1}
S.~Paul and G.~Tian.
\newblock {C}{M} {S}tability and the {G}eneralised {F}utaki {I}nvariant {I}.
\newblock {\em arXiv}, math.AG/0605278, 2006.

\bibitem{cms2}
S.~Paul and G.~Tian.
\newblock {C}{M} {S}tability and the {G}eneralised {F}utaki {I}nvariant {II}.
\newblock {\em arXiv}, math.AG/0606505, 2006.

\bibitem{gacms}
Sean~Timothy Paul.
\newblock Geometric analysis of {C}how {M}umford stability.
\newblock {\em Adv. Math.}, 182(2):333--356, 2004.

\bibitem{ruan}
Wei-Dong Ruan.
\newblock Canonical coordinates and {B}ergmann [{B}ergman] metrics.
\newblock {\em Comm. Anal. Geom.}, 6(3):589--631, 1998.

\bibitem{tianberg}
Gang Tian.
\newblock On a set of polarized {K}\"ahler metrics on algebraic manifolds.
\newblock {\em J. Differential Geom.}, 32(1):99--130, 1990.

\bibitem{kenhyp}
Gang Tian.
\newblock The {$K$}-energy on hypersurfaces and stability.
\newblock {\em Comm. Anal. Geom.}, 2(2):239--265, 1994.

\bibitem{psc}
Gang Tian.
\newblock K\"ahler-{E}instein metrics with positive scalar curvature.
\newblock {\em Invent. Math.}, 130(1):1--37, 1997.

\bibitem{eckart}
Eckart Viehweg.
\newblock {\em Quasi-projective moduli for polarized manifolds}, volume~30 of
  {\em Ergebnisse der Mathematik und ihrer Grenzgebiete (3) [Results in
  Mathematics and Related Areas (3)]}.
\newblock Springer-Verlag, Berlin, 1995.

\bibitem{zel}
Steve Zelditch.
\newblock Szeg{\H o} kernels and a theorem of {T}ian.
\newblock {\em Internat. Math. Res. Notices}, (6):317--331, 1998.

\bibitem{zhang}
Shouwu Zhang.
\newblock Heights and reductions of semi-stable varieties.
\newblock {\em Compositio Math.}, 104(1):77--105, 1996.

\end{thebibliography}
\end{document}